# Construction of analysis-suitable $G^1$ planar multi-patch parameterizations

Mario Kapl[a,b,*], Giancarlo Sangalli[b,c], Thomas Takacs[d]

[a]*Johann Radon Institute for Computational and Applied Mathematics, Austrian Academy of Sciences, Linz, Austria*
[b]*Dipartimento di Matematica "F.Casorati", Università degli Studi di Pavia, Italy*
[c]*Istituto di Matematica Applicata e Tecnologie Informatiche "E.Magenes" (CNR), Italy*
[d]*Institute of Applied Geometry, Johannes Kepler University Linz, Austria*

**Abstract**

Isogeometric analysis allows to define shape functions of global $C^1$ continuity (or of higher continuity) over multi-patch geometries. The construction of such $C^1$-smooth isogeometric functions is a non-trivial task and requires particular multi-patch parameterizations, so-called analysis-suitable $G^1$ (in short, AS-$G^1$) parameterizations, to ensure that the resulting $C^1$ isogeometric spaces possess optimal approximation properties, cf. [7]. In this work, we show through examples that it is possible to construct AS-$G^1$ multi-patch parameterizations of planar domains, given their boundary. More precisely, given a generic multi-patch geometry, we generate an AS-$G^1$ multi-patch parameterization possessing the same boundary, the same vertices and the same first derivatives at the vertices, and which is as close as possible to this initial geometry. Our algorithm is based on a quadratic optimization problem with linear side constraints. Numerical tests also confirm that $C^1$ isogeometric spaces over AS-$G^1$ multi-patch parameterized domains converge optimally under mesh refinement, while for generic parameterizations the convergence order is severely reduced.

## 1. Introduction

The concept of isogeometric analysis (IGA), cf. [2, 8, 14], provides the possibility to define smooth isogeometric spaces over multi-patch domains. The core idea is to use the same spline function space for representing the multi-patch geometry and for describing the solution of a partial differential equation. The high smoothness of the isogeometric spaces allows to solve high order partial differential equations directly via their weak form using a standard Galerkin discretization, see e.g. [1, 18, 19, 34].

---

[*]Corresponding author
*Email addresses:* mario.kapl@ricam.oeaw.ac.at (Mario Kapl), giancarlo.sangalli@unipv.it (Giancarlo Sangalli), thomas.takacs@jku.at (Thomas Takacs)

One important task in the framework of IGA is the construction of suitable multi-patch parameterizations for given multi-patch domains. In this paper, we deal with planar domains only. In this case, one possibility is to specify first the parameterizations of the boundaries of the single patches and to ensure at the same time that common interfaces have the same parameterizations. Then, the single patches are parameterized separately by fixing their boundaries. Several existing methods allow such a parameterization of single patches, e.g. [9, 10, 11, 13, 28]. If only the boundary of a multi-patch domain is given, the recent method [5] is able to generate a suitable multi-patch parameterization. The construction of globally $C^0$ smooth isogeometric functions over these $C^0$ multi-patch parameterizations is well understood and provides spaces with optimal approximation properties, see e.g. [2, 6, 31, 33]. But generic $C^0$ multi-patch parameterizations do not allow in general to generate $C^1$ isogeometric spaces with optimal approximation power, cf. [7].

Global $C^1$ continuity of isogeometric functions is equivalent, by definition, to the so-called geometric continuity $G^1$ of the graph of the function. This has been exploited in previous works cf. [7, 12, 20]. However, two different approaches, based on different types of multi-patch parameterizations, have been adopted. The first one requires multi-patch parameterizations which are $C^1$ almost everywhere except at extraordinary vertices (e.g. [26, 27, 32]), the second strategy generates $C^1$ isogeometric spaces over multi-patch parameterizations, which are only $C^0$ at the patch interfaces (e.g. [3, 4, 7, 16, 17, 20, 24, 25]). Similarly, [35] presents a construction of $C^k$ isogeometric spaces over patches having a polar layout, where they present details for $k \leq 2$.

In this paper, we are interested in the construction of $C^0$ multi-patch parameterizations which can be used when following the second strategy. It was recently shown in [7] that such multi-patch parameterizations have to satisfy specific constraints along the interfaces, derived from the geometric continuity condition, to guarantee that the resulting $C^1$ isogeometric spaces possess optimal approximation properties. Multi-patch parameterizations fulfilling exactly these constraints have been called analysis-suitable $G^1$ (in short, AS-$G^1$) multi-patch parameterizations in [7]. The class of AS-$G^1$ multi-patch parameterizations includes the subclass of bilinear multi-patch parameterizations (cf. [3, 16, 20, 24]) but the class of AS-$G^1$ multi-patch parameterization is wider than this subclass [7, 20].

However, already for generic biquadratic multi-patch parameterizations we obtain in general $C^1$ isogeometric spaces with dramatically reduced approximation properties. This is demonstrated by a first numerical example, presented in Fig. 1, which shows the importance of the use of AS-$G^1$ parameterizations. $L^2$ approximation is tested on a generic (initial) biquadratic polynomial $C^0$ multi-patch parameterization and on a bicubic spline AS-$G^1$ multi-patch parameterization possessing the same boundary. Note that the two parameterizations look indistinguishable in the plot. Nevertheless, in case of the generic parameterization the numerical solution does not even converge towards the exact one in the $L^\infty$-norm and converges with minimal order $\mathcal{O}(h^{1/2})$ in the $L^2$-norm. Whereas, in case of the AS-$G^1$ parameterization, the exact solution converges with an optimal rate of order $\mathcal{O}(h^4)$ in the $L^2$-norm.

The main goal of this paper is to show that it is possible to construct an AS-$G^1$ multi-patch parameterization for a given planar multi-patch domain by interpolating its boundary and by approximating its inner edges as good as possible. Though we do not develop a theoretical proof of this claim, we investigate numerically the following problem:

> Given a planar multi-patch domain $\widetilde{\Omega} = \cup_{\ell=1}^{P} \widetilde{\Omega}^{(\ell)}$ with an associated multi-patch parameterization $\widetilde{\boldsymbol{F}}$ consisting of single, regular parameterizations $\widetilde{\boldsymbol{F}}^{(\ell)}$, $\ell \in \{1, \ldots, P\}$. Find a multi-patch parameterization $\boldsymbol{F}$ consisting of single, regular parameterizations $\boldsymbol{F}^{(\ell)}$, $\ell \in \{1, \ldots, P\}$, which is AS-$G^1$, interpolates the boundary of $\widetilde{\boldsymbol{F}}$, interpolates the vertices and the first derivatives at the vertices of the single patches $\widetilde{\boldsymbol{F}}^{(\ell)}$, and approximates $\widetilde{\boldsymbol{F}}$ as good as possible. The resulting multi-patch parameterization $\boldsymbol{F}$ and resulting domain $\Omega = \cup_{\ell=1}^{P} \Omega^{(\ell)}$ represent the boundary of the initial multi-patch domain $\widetilde{\Omega}$ exactly and approximate the inner edges.

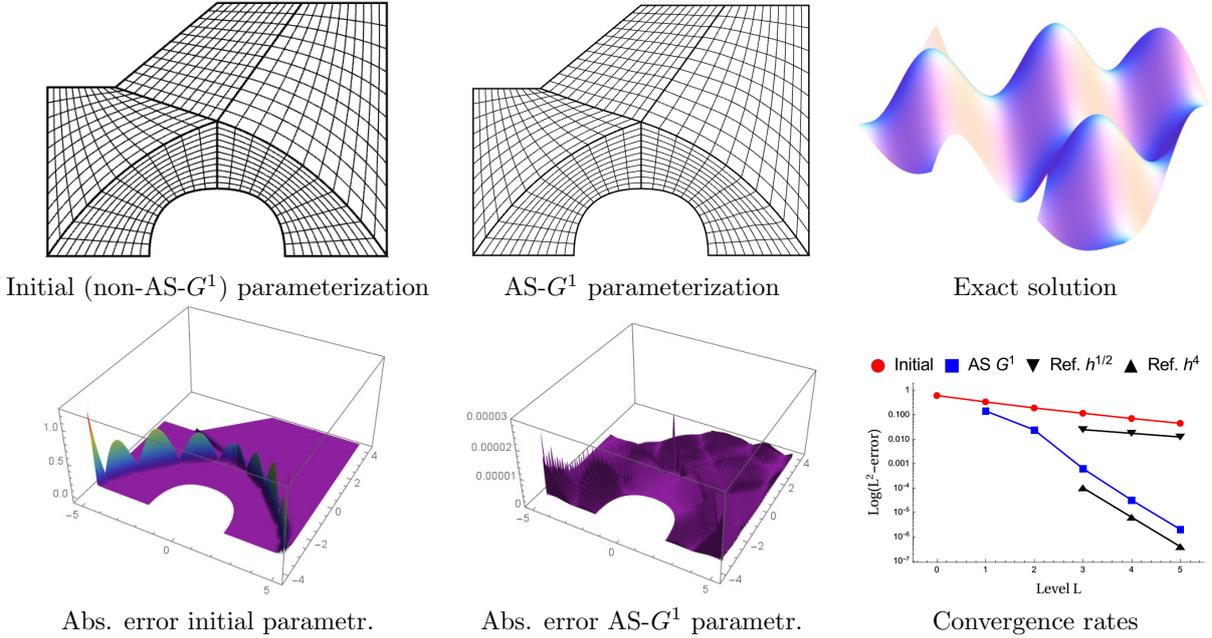

Figure 1: $L^2$ approximation of a given function (exact solution) over a non-AS-$G^1$ (initial) parameterization and over an AS-$G^1$ parameterization using $C^1$ isogeometric spline spaces of bidegree $(3,3)$. The numerical results indicate in case of the AS-$G^1$ parameterization optimal rates of convergence, and in case of the non-AS-$G^1$ parameterization no convergence in the $L^\infty$-norm and minimal rates of convergence in the $L^2$-norm. For both cases $C^1$ isogeometric spline spaces of bidegree $(3,3)$ are used for the approximation. See Section 4 for more detail about this example.

Observe that the problem above differs from the classical problem of a smooth interpolation of a mesh of curves, as it is described and studied in [30]. Indeed, we do not assume that the parametrization of the patch boundaries (i.e., the mesh of curves) is fixed. This condition would be too strong when coupled with the AS-$G^1$ constraint, as can be seen in numerical experiments. In fact, even the interpolation conditions at the vertices are not necessary in the configuration we consider. However, when omitting these vertex constraints, the problem becomes highly nonlinear and practically unsolveable. Certain relaxations of the constraints may be possible but beyond the scope of this paper.

In [16, 20] a method was proposed to construct mapped piecewise bilinear parameterizations, which are instances of AS-$G^1$ multi-patch parameterizations beyond bilinear multi-patch parameterizations. The two coordinate functions of such non-bilinear AS-$G^1$ parameterizations are constructed via least-squares fitting by means of $C^1$ isogeometric spaces over preselected bilinear domains. In contrast to our current approach the constructed AS-$G^1$ multi-patch parameterization does not interpolate the boundary, the vertices and the first derivatives at the vertices of a given multi-patch parameterization. Note that for certain domains no AS-$G^1$ parameterizations can be obtained with this approach, such as domains possessing a smooth boundary. Moreover, the selection (and influence) of the underlying bilinear parameterization for the $C^1$ isogeometric space was not investigated.

Different to our problem, in [4] a method was described to approximate a $G^0$ multi-patch surface by a $G^1$ multi-patch surface, where the gluing functions between the single patches are of low degree. The construction of the underlying $G^1$-smooth splines is based on the concept of syzygies. But so far, no isogeometric simulation has been performed to verify that the corresponding $C^1$ isogeometric spaces over these surfaces possess optimal approximation power.

To solve the issue stated above we formulate a quadratic minimization problem with linear side constraints. Therefore, the proposed strategy is simple and requires to solve only a system of linear equations.

We want to point out an important corollary of this work: Several examples of constructed AS-$G^1$ multi-patch parameterizations give clear evidence of the flexibility of AS-$G^1$ parameterization over planar multi-patch domains.

Furthermore, we perform numerical tests to analyze the order of convergence of $C^1$ isogeometric spaces with generic parameterization and AS-$G^1$ parameterizations (of the same planar domains). These test confirm the claims of [7], that is, optimal convergence is precluded in general on arbitrary parameterizations and is guaranteed on AS-$G^1$ parameterizations. These new tests are also more accurate that the ones in [7] since the $C^1$ isogeometric space is obtained by a symbolic null-space computation (using an approach based on minimal determining sets) rather than a numerical null-space computation.

The remainder of this paper is organized as follows. In Section 2, we describe some basic definitions and notations, and recall the concept of AS-$G^1$ multi-patch parameterizations. Section 3 presents the simple, linear method for constructing AS-$G^1$ multi-patch parameterizations. Thereby, the resulting AS-$G^1$ multi-patch geometry interpolates the boundary and the vertices and the first derivatives at the vertices of an initial multi-patch geometry by being at the same time as close as possible to this initial geometry. In Section 4, we construct several AS-$G^1$ multi-patch parameterizations for given multi-patch geometries to show the potential of our algorithm. Furthermore, we perform $L^2$ approximation on some of the constructed AS-$G^1$ parameterizations and their associated initial (non-AS-$G^1$) paramterization to verify that only AS-$G^1$ parameterizations guarantee $C^1$ isogeometric spaces with optimal approximation properties. Finally, we conclude the paper in Section 5.

## 2. Analysis-suitable $G^1$ multi-patch parameterizations

We follow a notation similar to the one introduced in [7], which we recall here. Let $\omega$ be the interval $[0,1]$ or the unit square $[0,1]^2$. The space $\mathcal{S}_k^{p,r}(\omega)$ is the (tensor-product) spline space of degree $p$ (in each direction) and continuity $C^r$ at all inner knots, which is defined on $\omega$ by choosing (in each direction) $k$ uniform inner knots of multiplicity $p-r$, where $k \in \mathbb{Z}_0^+$ and $0 \leq r \leq p$. Let $n = p + k(p-r)$. The spline spaces $\mathcal{S}_k^{p,r}([0,1])$ and $\mathcal{S}_k^{p,r}([0,1]^2)$ are spanned by the B-splines $N_i^{p,r}$, $i = 0, \ldots, n$, and the tensor-product B-splines $N_{i,j}^{p,r} = N_i^{p,r} N_j^{p,r}$, $i,j = 0, \ldots, n$, respectively. Moreover, we denote by $\mathcal{P}^p(w)$ the space of (tensor-product) polynomials of degree $p$ (in each direction) on $\omega$.

Consider a planar multi-patch domain $\Omega \subset \mathbb{R}^2$ which consists of $P \in \mathbb{Z}^+$ patches $\Omega^{(\ell)}$, i.e.

$$\Omega = \bigcup_{\ell=1}^{P} \Omega^{(\ell)},$$

where any two patches overlap only at the boundary

$$\Omega^{(\ell)} \cap \Omega^{(\ell')} \subset \partial \Omega^{(\ell)}$$

for all $\ell \neq \ell'$. Each patch $\Omega^{(\ell)}$, $\ell \in \{1, \ldots, P\}$, is the image $\boldsymbol{F}^{(\ell)}([0,1]^2)$ of a bijective and regular parameterization (also called *geometry mapping*)

$$\boldsymbol{F}^{(\ell)} : [0,1]^2 \to \Omega^{(\ell)}, \quad \boldsymbol{F}^{(\ell)} \in \mathcal{S}_k^{p,r}([0,1]^2) \times \mathcal{S}_k^{p,r}([0,1]^2),$$

with the spline representation

$$\boldsymbol{F}^{(\ell)}(u,v) = \sum_{i=0}^{n} \sum_{j=0}^{n} \boldsymbol{d}_{i,j}^{(\ell)} N_{i,j}^{p,r}(u,v), \quad \boldsymbol{d}_{i,j}^{(\ell)} \in \mathbb{R}^2. \tag{1}$$

We denote by $\boldsymbol{F}$ the multi-patch parameterization consisting of the parameterizations $\boldsymbol{F}^{(\ell)}$, $\ell \in \{1, \ldots, P\}$. In addition, we denote by $\Gamma^{(\ell,\ell')}$ the intersection of the two patches $\Omega^{(\ell)}$ and $\Omega^{(\ell')}$, $\ell, \ell' \in \{1, \ldots, P\}$ with $\ell \neq \ell'$, i.e.

$$\Gamma^{(\ell,\ell')} = \Omega^{(\ell)} \cap \Omega^{(\ell')}.$$

For simplicity, we assume that the intersection $\Gamma^{(\ell,\ell')}$ is either empty, a common vertex or a whole common boundary curve. Furthermore, we call two different patches $\Omega^{(\ell)}$ and $\Omega^{(\ell')}$ (or equivalently the two corresponding parameterizations $\boldsymbol{F}^{(\ell)}$ and $\boldsymbol{F}^{(\ell')}$), that share a whole common boundary curve, neighboring patches. We also assume that for each pair of neighboring patches $\Omega^{(\ell)}$ and $\Omega^{(\ell')}$ we (re)parameterize (by abuse of notation) $\boldsymbol{F}^{(\ell)}$ and $\boldsymbol{F}^{(\ell')}$ in such a way that

$$\boldsymbol{F}^{(\ell)}(0,v) = \boldsymbol{F}^{(\ell')}(0,v), \quad v \in [0,1], \tag{2}$$

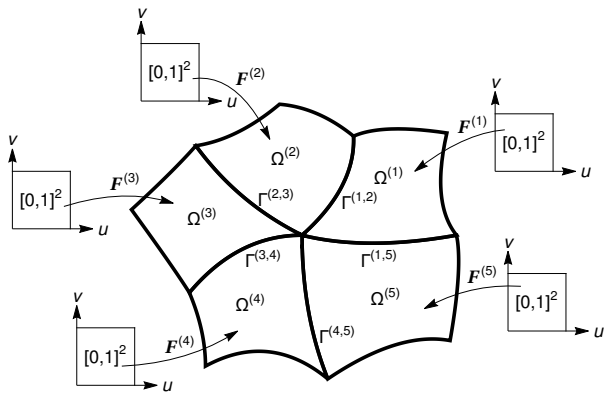

Figure 2: Example of the general setting of a multi-patch domain $\Omega$.

and denote the parameterization of the common curve (2) by $\boldsymbol{F}^{(\ell,\ell')}: [0,1] \to \mathbb{R}^2$. An example of a multi-patch domain $\Omega$ is visualized in Fig. 2.

The space $\mathcal{V}^1$ of $C^1$-smooth isogeometric functions on the multi-patch domain $\Omega$ with respect to the multi-patch parameterization $\boldsymbol{F}$ is given by

$$\mathcal{V}^1 = \{\phi \in C^1(\Omega) \mid \phi \circ \boldsymbol{F}^{(\ell)} \in \mathcal{S}_k^{p,r}([0,1]^2), \ \ell = 1, \ldots, P\},$$

see e.g. [7, 12, 16, 20]. Such spaces are of great interest and importance since they provide a possibility to solve numerically fourth order partial differential equations such as the biharmonic equation, see e.g. [16, 20]. To obtain $C^1$ isogeometric spaces with optimal approximation properties, particular multi-patch geometries, so-called analysis-suitable $G^1$ multi-patch geometries, are needed, cf. [7]. Let us recall this specific class of multi-patch parameterizations.

**Definition 1** (Analysis-suitable $G^1$ two-patch parameterization, cf. [7]). Two neighboring parameterizations $\boldsymbol{F}^{(\ell)}$ and $\boldsymbol{F}^{(\ell')}$ are called *analysis-suitable $G^1$* (in short, AS-$G^1$) at the interface $\Gamma^{(\ell,\ell')}$ if there exist $\alpha^{(\ell)}, \alpha^{(\ell')}, \beta^{(\ell)}, \beta^{(\ell')} \in \mathcal{P}^1([0,1])$ such that for all $v \in [0,1]$

$$\alpha^{(\ell)}(v)\alpha^{(\ell')}(v) < 0 \tag{3}$$

and

$$\alpha^{(\ell')}(v)D_u\boldsymbol{F}^{(\ell)}(0,v) - \alpha^{(\ell)}(v)D_u\boldsymbol{F}^{(\ell')}(0,v) + \beta(v)D_v\boldsymbol{F}^{(\ell,\ell')}(v) = \boldsymbol{0}, \tag{4}$$

with

$$\beta(v) = \alpha^{(\ell)}(v)\beta^{(\ell')}(v) - \alpha^{(\ell')}(v)\beta^{(\ell)}(v). \tag{5}$$

The functions $\alpha^{(\ell)}, \alpha^{(\ell')}, \beta^{(\ell)}, \beta^{(\ell')}$, and $\beta$ itself, are denoted as *gluing data* for the interface $\Gamma^{(\ell,\ell')}$.

Note that this is again an abuse of notation. Indeed, we require the existence of different $\alpha^{(\ell)}, \alpha^{(\ell')}, \ldots, \beta$ for every interface $\Gamma^{(\ell,\ell')}$. Whenever the AS-$G^1$ conditions are discussed, we will use the local notation above to keep the presentation simple.

The condition of AS-$G^1$ continuity is more restrictive than the classical $G^1$ continuity[1], where the gluing data are not required to be affine linear functions.

Equation (4) already uniquely determines $\alpha^{(\ell)}$, $\alpha^{(\ell')}$ and $\beta$ up to a common function $\gamma : [0,1] \to \mathbb{R}$ (with $\gamma(v) \neq 0$) by

$$\alpha^{(\ell)} = \gamma \bar{\alpha}^{(\ell)}, \alpha^{(\ell')} = \gamma \bar{\alpha}^{(\ell')} \text{ and } \beta = \gamma \bar{\beta}, \tag{6}$$

with $\bar{\alpha}^{(\ell)}, \bar{\alpha}^{(\ell')} \in \mathcal{S}_k^{2p-1,r-1}([0,1])$ and $\bar{\beta} \in \mathcal{S}_k^{2p,r}([0,1])$, where

$$\bar{\alpha}^{(\ell)}(v) = \det(D_u \boldsymbol{F}^{(\ell)}(0,v), D_v \boldsymbol{F}^{(\ell,\ell')}(v)), \tag{7}$$

$$\bar{\alpha}^{(\ell')}(v) = \det(D_u \boldsymbol{F}^{(\ell')}(0,v), D_v \boldsymbol{F}^{(\ell,\ell')}(v)), \tag{8}$$

and

$$\bar{\beta}(v) = \det(D_u \boldsymbol{F}^{(\ell)}(0,v), D_u \boldsymbol{F}^{(\ell')}(0,v)). \tag{9}$$

In general, there exists a one-parametric family of functions $\beta^{(\ell)}$ and $\beta^{(\ell')}$ such that (5) holds, see e.g. [7, 29].

**Definition 2** (Analysis-suitable $G^1$ multi-patch parameterization, cf. [7]). A multi-patch parameterization $\boldsymbol{F}$ is called *analysis-suitable* $G^1$ (in short, AS-$G^1$) if any two neighboring parameterizations $\boldsymbol{F}^{(\ell)}$ and $\boldsymbol{F}^{(\ell')}$ are AS-$G^1$-continuous at the corresponding interfaces $\Gamma^{(\ell,\ell')}$.

## 3. Construction of analysis-suitable $G^1$ multi-patch parameterizations

We present an algorithm for the construction of an AS-$G^1$ multi-patch parameterization from an initial multi-patch geometry.

### 3.1. Problem statement and outline

We consider the following problem: Given a planar multi-patch domain

$$\widetilde{\Omega} = \bigcup_{\ell=1}^{P} \widetilde{\Omega}^{(\ell)}$$

with an associated multi-patch parameterization $\widetilde{\boldsymbol{F}}$ consisting of regular geometry mappings $\widetilde{\boldsymbol{F}}^{(\ell)} \in \mathcal{S}_{\tilde{k}}^{\tilde{p},\tilde{r}}([0,1]^2) \times \mathcal{S}_{\tilde{k}}^{\tilde{p},\tilde{r}}([0,1]^2)$ with $\tilde{p} \geq 1$, $1 \leq \tilde{r} \leq \tilde{p}$, and $\tilde{k} \geq 0$. As before, we assume that for each pair of neighboring patches $\widetilde{\Omega}^{(\ell)}$ and $\widetilde{\Omega}^{(\ell')}$ we (re)parameterize (by abuse of notation) $\widetilde{\boldsymbol{F}}^{(\ell)}$ and $\widetilde{\boldsymbol{F}}^{(\ell')}$ in such a way that

$$\widetilde{\boldsymbol{F}}^{(\ell)}(0,v) = \widetilde{\boldsymbol{F}}^{(\ell')}(0,v), \quad v \in [0,1]. \tag{10}$$

---

[1] The definition of $G^1$ continuity of planar parametrizations is of little interest in geometric design; in cases of surfaces, see for example [29, 30].

The goal is to find a multi-patch parameterization $\boldsymbol{F}$ consisting of geometry mappings $\boldsymbol{F}^{(\ell)} \in \mathcal{S}_k^{p,r}([0,1]^2) \times \mathcal{S}_k^{p,r}([0,1]^2)$ with $\mathcal{S}_k^{p,r}([0,1]^2) \supseteq \mathcal{S}_{\tilde{k}}^{\tilde{p},\tilde{r}}([0,1]^2)$, and $1 \leq r \leq p-2$, which describes the planar multi-patch domain

$$\Omega = \bigcup_{\ell=1}^{P} \Omega^{(\ell)},$$

such that $\boldsymbol{F}$ approximates $\widetilde{\boldsymbol{F}}$ as good as possible, i.e.

$$\boldsymbol{F} \approx \widetilde{\boldsymbol{F}},$$

and the following conditions are satisfied:

(a) $\boldsymbol{F}$ is AS-$G^1$,

(b) the boundaries of the two multi-patch domains $\Omega$ and $\widetilde{\Omega}$ coincide, i.e.

$$\partial \Omega = \partial \widetilde{\Omega}, \tag{11}$$

and their parameterizations by $\boldsymbol{F}$ and $\widetilde{\boldsymbol{F}}$ coincide,

(c) for all $\ell \in \{1, \ldots, P\}$ the vertices and the first derivatives at the vertices of the parameterizations $\boldsymbol{F}^{(\ell)}$ and $\widetilde{\boldsymbol{F}}^{(\ell)}$ coincide, i.e.

$$\boldsymbol{F}^{(\ell)}(u_0, v_0) = \widetilde{\boldsymbol{F}}^{(\ell)}(u_0, v_0), \quad u_0, v_0 \in \{0, 1\}, \tag{12}$$

and

$$D_u \boldsymbol{F}^{(\ell)}(u_0, v_0) = D_u \widetilde{\boldsymbol{F}}^{(\ell)}(u_0, v_0), \quad u_0, v_0 \in \{0, 1\}, \tag{13}$$

$$D_v \boldsymbol{F}^{(\ell)}(u_0, v_0) = D_v \widetilde{\boldsymbol{F}}^{(\ell)}(u_0, v_0), \quad u_0, v_0 \in \{0, 1\}, \tag{14}$$

respectively, and

(d) for all $\ell \in \{1, \ldots, P\}$ the parameterization $\boldsymbol{F}^{(\ell)}$ is regular.

Below, we will refer to these conditions as conditions (a)-(d).

**Remark 1.** Note that since we consider planar domains conditions (a)-(d) guarantee that $\Omega = \widetilde{\Omega}$. When extending the construction to surfaces in $\mathbb{R}^3$, cf. [7], this is not satisfied anymore. The extension to surface domains is straight-forward, but omitted here to keep the presentation short. A study on the approximation quality in that case is of vital interest for future research.

Recall that each parameterization $\boldsymbol{F}^{(\ell)} \in \mathcal{S}_k^{p,r}([0,1]^2) \times \mathcal{S}_k^{p,r}([0,1]^2)$, $\ell \in \{1, \ldots, P\}$, has a B-spline representation (1) with control points $\boldsymbol{d}_{i,j}^{(\ell)} \in \mathbb{R}^2$. In order to simplify the notation, we denote by $\boldsymbol{d}$ the matrix of all control points $\boldsymbol{d}_{i,j}^{(\ell)}$ from all parameterizations $\boldsymbol{F}^{(\ell)}$, $\ell \in \{1, \ldots, P\}$, i.e.

$$\boldsymbol{d} = (\boldsymbol{d}_{i,j}^{(\ell)})_{\ell=1,\ldots,P; i,j=0,\ldots,n}.$$

We generate the desired planar multi-patch parameterization $\boldsymbol{F}$ by minimizing the following objective function

$$\mathcal{F}_2(\boldsymbol{d}) + \lambda_L \mathcal{F}_L(\boldsymbol{d}) + \lambda_U \mathcal{F}_U(\boldsymbol{d}) \to \min_{\boldsymbol{d}} \tag{15}$$

with respect to specific linear constraints $\mathcal{A}_G$, $\mathcal{A}_B$ and $\mathcal{A}_V$, where $\lambda_L$ and $\lambda_U$ are user-defined non-negative weights, and the quadratic functionals $\mathcal{F}_2$, $\mathcal{F}_L$ and $\mathcal{F}_U$ are given by

$$\mathcal{F}_2(\boldsymbol{d}) = \sum_{\ell=1}^{P} \int_{[0,1]^2} \|\boldsymbol{F}^{(\ell)} - \widetilde{\boldsymbol{F}}^{(\ell)}\|^2 \, \mathrm{d}u \, \mathrm{d}v,$$

$$\mathcal{F}_L(\boldsymbol{d}) = \sum_{\ell=1}^{P} \int_{[0,1]^2} \left( \|D_u \boldsymbol{F}^{(\ell)}\|^2 + \|D_v \boldsymbol{F}^{(\ell)}\|^2 \right) \mathrm{d}u \, \mathrm{d}v,$$

and

$$\mathcal{F}_U(\boldsymbol{d}) = \sum_{\ell=1}^{P} \int_{[0,1]^2} \left( \|D_{uu} \boldsymbol{F}^{(\ell)}\|^2 + 2\|D_{uv} \boldsymbol{F}^{(\ell)}\|^2 + \|D_{vv} \boldsymbol{F}^{(\ell)}\|^2 \right) \mathrm{d}u \, \mathrm{d}v,$$

respectively. The linear constraints $\mathcal{A}_G$, $\mathcal{A}_B$ and $\mathcal{A}_V$ are used to satisfy condition (a), (b) and (c), respectively. The derivation of these linear constraints will be explained in the following subsection. The functional $\mathcal{F}_2$ guarantees that the resulting multi-patch parameterization $\boldsymbol{F}$ is close to the initial multi-patch parameterization $\widetilde{\boldsymbol{F}}$. The functionals $\mathcal{F}_L$ and $\mathcal{F}_U$ (the so-called parametric length functional and the uniformity functional, cf. [5, 9]) are used to obtain parameterizations $\boldsymbol{F}$ of good quality, and are controlled by the non-negative weights $\lambda_L$ and $\lambda_U$, respectively.

The regularity condition (d) is nonlinear and cannot be guaranteed directly by our method. However, it can be validated after optimization. We conjecture that a sufficiently refined spline space always allows for regular parameterizations $\boldsymbol{F}^{(\ell)}$. Note that all the examples we present here satisfy condition (d).

The concept of least-squares fitting is a well known tool to generate (surface) parameterizations, e.g. [36, 37]. The proposed linear direct method was recently used in a similar way in [5, 9] to parameterize (multi-patch) domains. The novelty of our approach is the use of specific linear constraints to ensure that the resulting multi-patch parameterization $\boldsymbol{F}$ is AS-$G^1$.

*3.2. The linear constraints*

We describe the linear constraints $\mathcal{A}_G$, $\mathcal{A}_B$ and $\mathcal{A}_V$, which are used to fulfill condition (a), (b) and (c), respectively. In order to simplify the computation, we assume that the parameterizations $\boldsymbol{F}^{(\ell)}$ and $\widetilde{\boldsymbol{F}}^{(\ell)}$, $\ell \in \{1, \ldots, P\}$, are given by B-spline representations (1) for the same space $\mathcal{S}_k^{p,r}([0,1]^2)$, more precisely,

$$\boldsymbol{F}^{(\ell)}(u,v) = \sum_{i=0}^{n} \sum_{j=0}^{n} \boldsymbol{d}_{i,j}^{(\ell)} N_{i,j}^{p,r}(u,v), \quad \boldsymbol{d}_{i,j}^{(\ell)} \in \mathbb{R}^2,$$

and
$$\widetilde{\boldsymbol{F}}^{(\ell)}(u,v) = \sum_{i=0}^{n}\sum_{j=0}^{n} \widetilde{\boldsymbol{d}}_{i,j}^{(\ell)} N_{i,j}^{p,r}(u,v), \quad \widetilde{\boldsymbol{d}}_{i,j}^{(\ell)} \in \mathbb{R}^2. \tag{16}$$

In case of $\mathcal{S}_{\widetilde{k}}^{\widetilde{p},\widetilde{r}}([0,1]^2) \subsetneq \mathcal{S}_{k}^{p,r}([0,1]^2)$, the B-spline representation (16) can be simply obtained from the given B-spline representation for the space $\mathcal{S}_{\widetilde{k}}^{\widetilde{p},\widetilde{r}}([0,1]^2)$ by a suitable raise of degree and/or insertion of knots.

*AS-$G^1$ constraints $\mathcal{A}_G$.* To ensure that the multi-patch parameterization is AS-$G^1$, compare condition (a), we have to satisfy for all interfaces $\Gamma^{(\ell,\ell')}$ that (2) holds, and that there exists gluing data $\alpha^{(\ell)}, \alpha^{(\ell')}, \beta^{(\ell)}, \beta^{(\ell')} \in \mathcal{P}^1([0,1])$ such that (3), (4) and (5) are fulfilled. By considering these constraints more closely, we can easily observe that we obtain highly non-linear ones. The reason is that on the one hand the gluing data $\alpha^{(\ell)}, \alpha^{(\ell')}, \beta^{(\ell)}$ and $\beta^{(\ell')}$ are not given, and on the other hand most of the control points $\boldsymbol{d}_{i,j}^{(\ell)}$ (except the ones fixed by means of constraints $\mathcal{A}_B$ and $\mathcal{A}_V$), which are affected by the AS-$G^1$ conditions (2), (3), (4) and (5), are unknown, too.

Therefore, we simplify these constraints to linear ones by first computing the linear gluing data $\alpha^{(\ell)}, \alpha^{(\ell')}, \beta^{(\ell)}$ and $\beta^{(\ell')}$ by means of the initial multi-patch parameterization $\widetilde{\boldsymbol{F}}$ and the relations (6)-(9). Consider an interface $\Gamma^{(\ell,\ell')}$. The aim is to construct linear functions
$$\alpha^{(\ell)}(v) = a_0^{(\ell)}(1-v) + a_1^{(\ell)}v, \quad a_0^{(\ell)}, a_1^{(\ell)} \in \mathbb{R},$$
$$\alpha^{(\ell')}(v) = a_0^{(\ell')}(1-v) + a_1^{(\ell')}v, \quad a_0^{(\ell')}, a_1^{(\ell')} \in \mathbb{R},$$
and
$$\beta^{(\ell)}(v) = b_0^{(\ell)}(1-v) + b_1^{(\ell)}v, \quad b_0^{(\ell)}, b_1^{(\ell)} \in \mathbb{R},$$
$$\beta^{(\ell')}(v) = b_0^{(\ell')}(1-v) + b_1^{(\ell')}v, \quad b_0^{(\ell')}, b_1^{(\ell')} \in \mathbb{R},$$
which fit well with the functions
$$\widetilde{\alpha}^{(\ell)}(v) = \det(D_u\widetilde{\boldsymbol{F}}^{(\ell)}(0,v), D_v\widetilde{\boldsymbol{F}}^{(\ell,\ell')}(v)),$$
$$\widetilde{\alpha}^{(\ell')}(v) = \det(D_u\widetilde{\boldsymbol{F}}^{(\ell')}(0,v), D_v\widetilde{\boldsymbol{F}}^{(\ell,\ell')}(v)),$$
and
$$\widetilde{\beta}(v) = \det(D_u\widetilde{\boldsymbol{F}}^{(\ell)}(0,v), D_u\widetilde{\boldsymbol{F}}^{(\ell')}(0,v)).$$
obtained from the initial parameterizations $\widetilde{\boldsymbol{F}}^{(\ell)}$ and $\widetilde{\boldsymbol{F}}^{(\ell')}$. More precisely, we generate the gluing data $\alpha^{(\ell)}$ and $\alpha^{(\ell')}$ by choosing their control points as
$$a_0^{(\ell)} = \widetilde{\alpha}^{(\ell)}(0), \; a_1^{(\ell)} = \widetilde{\alpha}^{(\ell)}(1), \; a_0^{(\ell')} = \widetilde{\alpha}^{(\ell')}(0), \; a_1^{(\ell')} = \widetilde{\alpha}^{(\ell')}(1),$$
and the functions $\beta^{(\ell)}$ and $\beta^{(\ell')}$ by minimizing the objective function
$$\begin{aligned}&\int_0^1 \|\widetilde{\beta} - \underbrace{(\alpha^{(\ell)}\beta^{(\ell')} - \alpha^{(\ell')}\beta^{(\ell)})}_{\beta}\|^2 \mathrm{d}v \\ &+ \lambda_\beta \left(\int_0^1 \|\beta^{(\ell)}\|^2 \mathrm{d}v + \int_0^1 \|\beta^{(\ell')}\|^2 \mathrm{d}v\right) \to \min_{(b_0^{(\ell)}, b_1^{(\ell)}, b_0^{(\ell')}, b_0^{(\ell')})}\end{aligned} \tag{17}$$

with respect to the linear constraints

$$\beta(0) = \widetilde{\beta}(0) \text{ and } \beta(1) = \widetilde{\beta}(1),$$

where $\lambda_\beta$ is a user-defined positive weight. The second and third term in the objective function (17) are needed to obtain a non-singular linear system, and the influence of these terms in the minimization process are controlled by the positive weight $\lambda_\beta$. The selection of the gluing data $\alpha^{(\ell)}$, $\alpha^{(\ell')}$, $\beta^{(\ell)}$ and $\beta^{(\ell')}$ guarantees that

$$\alpha^{(\ell)}(v_0) = \widetilde{\alpha}^{(\ell)}(v_0),\ \alpha^{(\ell')}(v_0) = \widetilde{\alpha}^{(\ell')}(v_0),\ \beta(v_0) = \widetilde{\beta}(v_0),$$

for $v_0 \in \{0,1\}$, which later allows to state the linear constraints $\mathcal{A}_B$ and $\mathcal{A}_V$. Note that this is equivalent to setting $\gamma(v) \equiv 1$ in (6). Our numerical tests suggest that a different choice has no significant influence on the fitting result. Furthermore, the choice of $\alpha^{(\ell)}$ and $\alpha^{(\ell')}$ ensures that (3) holds. Finally, for all interfaces $\Gamma^{(\ell,\ell')}$ the conditions (2), (4) and (5) with the selected gluing data $\alpha^{(\ell)}$, $\alpha^{(\ell')}$, $\beta^{(\ell)}$ and $\beta^{(\ell')}$ lead to desired linear constraints on the involved control points $\boldsymbol{d}_{i,j}^{(\ell)}$.

*Boundary constraints $\mathcal{A}_B$.* Let $\boldsymbol{b}_m : [0,1] \to [0,1]^2$, $m \in \{1, \ldots, 4\}$, be the functions given by

$$\boldsymbol{b}_m(t) = \begin{cases} (t,0) & \text{if } m = 1, \\ (t,1) & \text{if } m = 2, \\ (0,t) & \text{if } m = 3, \\ (1,t) & \text{if } m = 4. \end{cases}$$

To fulfill condition (b), i.e., that the multi-patch parameterizations $\boldsymbol{F}$ and $\widetilde{\boldsymbol{F}}$ coincide at the domain boundary, or equivalently that for each $\ell \in \{1, \ldots, P\}$, and each function $\boldsymbol{b}_m$, $m \in \{1, \ldots, 4\}$, satisfying

$$\{\widetilde{\boldsymbol{F}}^{(\ell)}(\boldsymbol{b}_i(t)),\ t \in [0,1]\} \subseteq \partial\widetilde{\Omega},$$

the two boundary curves $\boldsymbol{F}^{(\ell)}(\boldsymbol{b}_m(t))$ and $\widetilde{\boldsymbol{F}}^{(\ell)}(\boldsymbol{b}_m(t))$ are equal,

$$\boldsymbol{F}^{(\ell)}(\boldsymbol{b}_m(t)) = \widetilde{\boldsymbol{F}}^{(\ell)}(\boldsymbol{b}_m(t)), \quad t \in [0,1],$$

we require that

$$\begin{cases} \boldsymbol{d}_{i,0}^{(\ell)} = \widetilde{\boldsymbol{d}}_{i,0}^{(\ell)}, & i = 0, \ldots, n, \quad \text{if } m = 1, \\ \boldsymbol{d}_{i,n}^{(\ell)} = \widetilde{\boldsymbol{d}}_{i,n}^{(\ell)}, & i = 0, \ldots, n, \quad \text{if } m = 2, \\ \boldsymbol{d}_{0,j}^{(\ell)} = \widetilde{\boldsymbol{d}}_{0,j}^{(\ell)}, & j = 0, \ldots, n, \quad \text{if } m = 3, \\ \boldsymbol{d}_{n,j}^{(\ell)} = \widetilde{\boldsymbol{d}}_{n,j}^{(\ell)}, & j = 0, \ldots, n, \quad \text{if } m = 4. \end{cases} \tag{18}$$

The resulting constraints are compatible with the constraints $\mathcal{A}_V$ and with the selection of the gluing data $\alpha^{(\ell)}$, $\alpha^{(\ell')}$, $\beta^{(\ell)}$ and $\beta^{(\ell')}$ for the single interfaces $\Gamma^{(\ell,\ell')}$ in constraints $\mathcal{A}_G$.

*Vertex constraints* $\mathcal{A}_V$. The aim is to fulfill for all $\ell \in \{1, \ldots, P\}$ the equations (12)-(14), compare condition (c). This can be simply obtained by defining the constraints

$$\boldsymbol{d}_{0,0}^{(\ell)} = \widetilde{\boldsymbol{d}}_{0,0}^{(\ell)}, \quad \boldsymbol{d}_{1,0}^{(\ell)} = \widetilde{\boldsymbol{d}}_{1,0}^{(\ell)}, \quad \boldsymbol{d}_{0,1}^{(\ell)} = \widetilde{\boldsymbol{d}}_{0,1}^{(\ell)}, \tag{19}$$

$$\boldsymbol{d}_{0,n}^{(\ell)} = \widetilde{\boldsymbol{d}}_{0,n}^{(\ell)}, \quad \boldsymbol{d}_{1,n}^{(\ell)} = \widetilde{\boldsymbol{d}}_{1,n}^{(\ell)}, \quad \boldsymbol{d}_{0,n-1}^{(\ell)} = \widetilde{\boldsymbol{d}}_{0,n-1}^{(\ell)}, \tag{20}$$

$$\boldsymbol{d}_{n,0}^{(\ell)} = \widetilde{\boldsymbol{d}}_{n,0}^{(\ell)}, \quad \boldsymbol{d}_{n-1,0}^{(\ell)} = \widetilde{\boldsymbol{d}}_{n-1,0}^{(\ell)}, \quad \boldsymbol{d}_{n,1}^{(\ell)} = \widetilde{\boldsymbol{d}}_{n,1}^{(\ell)}, \tag{21}$$

and

$$\boldsymbol{d}_{n,n}^{(\ell)} = \widetilde{\boldsymbol{d}}_{n,n}^{(\ell)}, \quad \boldsymbol{d}_{n-1,n}^{(\ell)} = \widetilde{\boldsymbol{d}}_{n-1,n}^{(\ell)}, \quad \boldsymbol{d}_{n,n-1}^{(\ell)} = \widetilde{\boldsymbol{d}}_{n,n-1}^{(\ell)}, \tag{22}$$

for $\ell \in \{1, \ldots, P\}$. These constraints are compatible with the constraints $\mathcal{A}_B$ and with the selection of the gluing data $\alpha^{(\ell)}$, $\alpha^{(\ell')}$, $\beta^{(\ell)}$ and $\beta^{(\ell')}$ for the single interfaces $\Gamma^{(\ell,\ell')}$ in constraints $\mathcal{A}_G$.

### 3.3. The construction method

We describe a symbolic algorithm for the construction of an AS-$G^1$ planar multi-patch parameterization $\boldsymbol{F}$ consisting of parameterizations $\boldsymbol{F}^{(\ell)}$. All involved integrals will be evaluated numerically (by means of the rectangle method) to speed up and to simplify the computation, which will be done completely symbolically.

The construction of the desired AS-$G^1$ multi-patch parameterization $\boldsymbol{F}$ works now as follows:

1. **Input:** We have given a planar multi-patch domain $\widetilde{\Omega} = \cup_{\ell=1}^{P} \widetilde{\Omega}^{(\ell)}$ with an associated multi-patch parameterization $\widetilde{\boldsymbol{F}}$ consisting of regular geometry mappings $\widetilde{\boldsymbol{F}}^{(\ell)} \in \mathcal{S}_{\widetilde{k}}^{\widetilde{p},\widetilde{r}}([0,1]^2) \times \mathcal{S}_{\widetilde{k}}^{\widetilde{p},\widetilde{r}}([0,1]^2)$, where $\widetilde{p} \geq 1$, $1 \leq \widetilde{r} \leq \widetilde{p}$, and $\widetilde{k} \geq 0$.

2. The space $\mathcal{S}_k^{p,r}([0,1]^2) \times \mathcal{S}_k^{p,r}([0,1]^2)$ for the single parameterizations $\boldsymbol{F}^{(\ell)}$ is selected by choosing the desired degree $p$, the regularity $r$ and the number $k$ of different inner knots. Thereby, it has to be satisfied that

$$\mathcal{S}_k^{p,r}([0,1]^2) \supseteq \mathcal{S}_{\widetilde{k}}^{\widetilde{p},\widetilde{r}}([0,1]^2),$$

   and $1 \leq r \leq p-2$.

3. The linear constraints $\mathcal{A}_B$ and $\mathcal{A}_V$ already determine the control points $\boldsymbol{d}_{i,j}^{(\ell)}$ given in (18) and in (19)-(22), respectively. Then, we compute the linear constraints $\mathcal{A}_G$ as described in Section 3.2. For this we select suitable positive weights $\lambda_\beta$ for the interfaces $\Gamma^{(\ell,\ell')}$. For simplicity, we choose for all interfaces the same weight $\lambda_\beta = \frac{1}{100}$ (larger values work fine too). The constraints $\mathcal{A}_G$ specify a system of linear equations, which is in general under-determined and does not have full rank. Let $q$ be this rank. We solve symbolically the system and obtain a solution, where $q$ coefficients of $\boldsymbol{d}$ are linearly described by other coefficients of $\boldsymbol{d}$ playing the role of free variables.

4. We select suitable non-negative weights $\lambda_L$ and $\lambda_U$, where a good choice of these weights for our examples (see Section 4) has been observed to lie within the interval $[0, \frac{1}{1000}]$. Recall that these weights control the magnitude of influence of the functionals $\mathcal{F}_L$ and $\mathcal{F}_U$ with respect to the functional $\mathcal{F}_2$. They have to be selected in such a way that on the one hand the resulting relative $L^2$ error is small and on the other hand the resulting parameterization is of good quality. Then, we compute the objective function (15) and minimize it with respect to the remaining unknown coefficients of $\boldsymbol{d}$. Again, we solve symbolically the resulting system of linear equations and obtain a unique solution for the remaining unknown coefficients of $\boldsymbol{d}$, which finally determine all control points $\boldsymbol{d}_{i,j}^{(\ell)}$ of the parameterizations $\boldsymbol{F}^{(\ell)}$, $\ell \in \{1, \ldots, P\}$.

5. **Output:** We obtain a multi-patch parameterization $\boldsymbol{F}$ with single parameterizations $\boldsymbol{F}^{(\ell)}$ and the corresponding multi-patch domain $\Omega = \cup_{\ell=1}^{P} \Omega_\ell$, such that conditions (a)-(d) are satisfied.

**Remark 2.** Our experiments clearly show that refining the space $\mathcal{S}_k^{p,r}([0,1]^2)$ (i.e., by increasing multiplicity of existing knots and/or by inserting new inner knots and/or by raising the degree) significantly improves the quality of the resulting multi-patch parameterization $\boldsymbol{F}$. Furthermore, the choice of the weights $\lambda_L$ and $\lambda_U$ affects the quality of the computed parameterization. Moreover, the parameterization may not be regular, i.e. condition (d) may not be satisfied. Also in that case, refining the function space improves the result.

**Remark 3.** Note that for our approach it is not necessary that an initial multi-patch parameterization $\widetilde{\boldsymbol{F}}$ is given. If only the multi-patch domain $\widetilde{\Omega}$ is given as a collection of patches $\Omega^{(\ell)}$, an initial parameterization $\widetilde{\boldsymbol{F}}$ can be computed as follows. We first specify and fix parameterizations for all interfaces and ensure that the patch parameterizations $\widetilde{\boldsymbol{F}}^{(\ell)}$ coincide with these interface parameterizations. Then we obtain an initial parameterization $\widetilde{\boldsymbol{F}}$ by minimizing the objective function

$$\lambda_L \mathcal{F}_L(\widetilde{\boldsymbol{d}}) + \lambda_U \mathcal{F}_U(\widetilde{\boldsymbol{d}}) \to \min_{\widetilde{\boldsymbol{d}}},$$

with respect to the fixed boundaries of the parameterizations $\widetilde{\boldsymbol{F}}^{(\ell)}$. For most standard CAD geometries, not even the patches $\Omega^{(\ell)}$ are given, but only the domain boundary $\partial \widetilde{\Omega}$. In that case the recent method [5] can be used to construct an initial multi-patch parameterization $\widetilde{\boldsymbol{F}}$.

**Remark 4.** Our tests have shown that the AS-$G^1$-constraints are highly sensitive to small perturbations. This has motivated the use of a symbolic (exact arithmetic) solver. Our current and further research will focus on the analytic construction of a local and well conditioned basis for the spaces that fulfill the AS-$G^1$ constraints $\mathcal{A}_G$, extending the results of [7, 16, 17, 20].

## 4. Examples

We first demonstrate the performance of our algorithm by constructing several examples of AS-$G^1$ planar multi-patch parameterizations. Then, we perform $L^2$ approximation on some of the resulting AS-$G^1$ parameterizations and on their underlying initial non-AS-$G^1$ parameterizations to verify on the one hand that AS-$G^1$ parameterizations allow to generate $C^1$ isogeometric spaces with optimal approximation properties and to show on the other hand that non-AS-$G^1$ parameterizations lead in general to $C^1$ isogeometric spaces with dramatically reduced approximation properties.

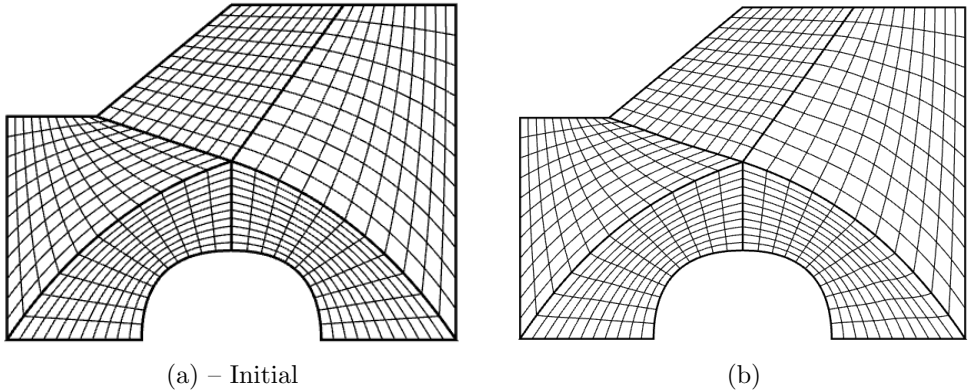

(a) – Initial                                                (b)

Figure 3: Example 1 – Part of car. Compare Table 1.

*4.1. Instances of AS-$G^1$ multi-patch parameterizations*

We use the algorithm presented in the previous section to generate AS-$G^1$ planar multi-patch parameterizations $\boldsymbol{F}$ from initial multi-patch parameterizations $\widetilde{\boldsymbol{F}}$. Thereby, the algorithm is implemented in Mathematica, which allows to perform all computations symbolically. The resulting AS-$G^1$ multi-patch parameterizations $\boldsymbol{F}$ are visualized in Figs. 3-6, and the parameters $p$, $r$ and $k$ for the selected spaces $\mathcal{S}_k^{p,r}([0,1]^2) \times \mathcal{S}_k^{p,r}([0,1]^2)$ as well as the choices of the non-negative weights $\lambda_L$ and $\lambda_U$ are specified in Table 1. In addition, the table reports the resulting relative $L^2$ errors and the resulting values for the functional $\mathcal{F}_2$, for the parametric length functional $\mathcal{F}_L$ and for the uniformity functional $\mathcal{F}_U$ in the objective function (15).

Below, we will consider the single examples in detail.

*Example 1 – Part of car.* The initial five-patch domain $\widetilde{\Omega}$ shown in Fig. 3(a) is inspired by [16, Fig. 12] and could represent a part of a car. Thereby, the single patches $\widetilde{\Omega}^{(\ell)}$ are represented by parameterizations $\widetilde{\boldsymbol{F}}^{(\ell)} \in \mathcal{P}^2([0,1]^2) \times \mathcal{P}^2([0,1]^2)$. The use of the space $\mathcal{S}_1^{3,1}([0,1]^2) \times \mathcal{S}_1^{3,1}([0,1]^2)$ for the parameterizations $\boldsymbol{F}^{(\ell)}$ lead to an AS-$G^1$ multi-patch parameterization $\boldsymbol{F}$ of good quality with small relative $L^2$ error, see Fig. 3(b) and Table 1. (Note that in Table 1 we refer to the space $\mathcal{P}^2([0,1]^2)$ as $\mathcal{S}_0^{2,\infty}([0,1]^2)$.) In contrast to the proposed method [16, Example 6] for this particular example we are able to interpolate the boundary, the vertices and the first derivatives at the vertices of the initial parameterization $\widetilde{\boldsymbol{F}}$.

*Example 2 – Puzzle piece.* This example is inspired by the domain shown in [5, Fig. 6], which is the quarter of a puzzle piece. We consider the initial multi-patch domain $\widetilde{\Omega}$ visualized in Fig. 4(a), which consists of 6 patches $\widetilde{\Omega}^{(\ell)}$. The single patches $\widetilde{\Omega}^{(\ell)}$ are parameterized by geometry mappings $\widetilde{\boldsymbol{F}}^{(\ell)} \in \mathcal{S}_2^{3,1}([0,1]^2) \times \mathcal{S}_2^{3,1}([0,1]^2)$. Fig. 4(b) shows the resulting AS-$G^1$ parameterizations $\boldsymbol{F}$ choosing the same space $\mathcal{S}_2^{3,1}([0,1]^2) \times \mathcal{S}_2^{3,1}([0,1]^2)$ for the single patches $\boldsymbol{F}^{(\ell)}$, and Table 1 reports the resulting relative $L^2$ error.

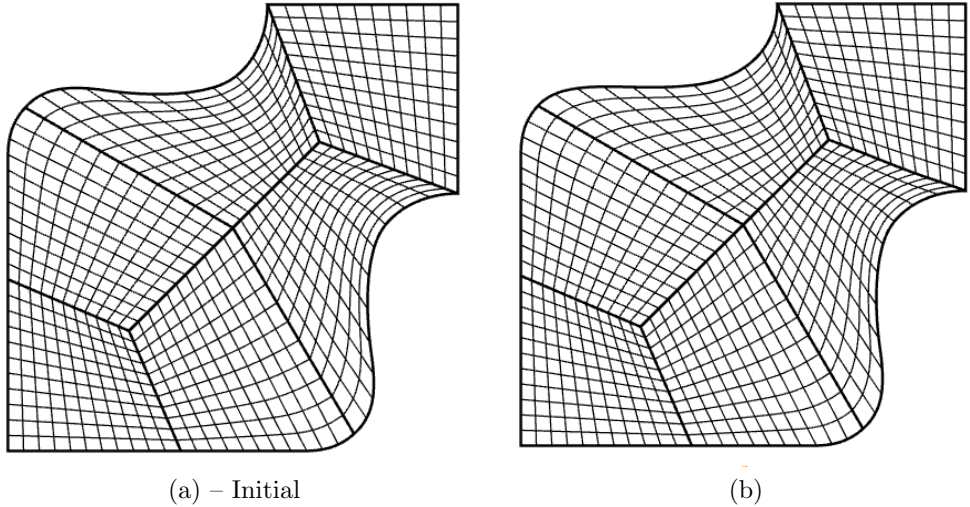

(a) – Initial                  (b)

Figure 4: Example 2 – Puzzle piece. Compare Table 1.

*Example 3 – Deformed circle.* We consider an initial multi-patch domain $\widetilde{\Omega}$ consisting of five patches $\widetilde{\Omega}^{(\ell)}$, where the single patches $\widetilde{\Omega}^{(\ell)}$ are parameterized by geometry mappings $\widetilde{\boldsymbol{F}}^{(\ell)} \in \mathcal{S}_2^{3,2}([0,1]^2) \times \mathcal{S}_2^{3,2}([0,1]^2)$, see Fig. 5(a). Thereby, the boundary $\partial\widetilde{\Omega}$ is parameterized in such a way that it is $C^1$ smooth. Note that in this example none of the edges of the single patches are straight lines. We compute AS-$G^1$ multi-patch parameterizations $\boldsymbol{F}$ for two different selections of the parameters $p$, $r$ and $k$ for the spaces $\mathcal{S}_k^{p,r}([0,1]^2)$, see Fig. 5(b)-(c) and Table 1, and for both cases the resulting AS-$G^1$ multi-patch geometries possess a parameterization of good quality with a small relative $L^2$ error.

In addition, we plot for the interface located in the blue encircled region the gluing data $\alpha^{(\ell)}$, $\alpha^{(\ell')}$ and $\beta$ of the initial parameterization (a) (in black) and of the obtained AS-$G^1$ parameterizations (b)-(c) (in red), see Fig 5(d)-(f). Note that each resulting AS-$G^1$ parameterization possess the same gluing data $\alpha^{(\ell)}$, $\alpha^{(\ell')}$ and $\beta$ independent of the used space $\mathcal{S}_k^{p,r}([0,1]^2) \times \mathcal{S}_k^{p,r}([0,1]^2)$ as a consequence of our construction method. One can verify that for the AS-$G^1$ parameterizations the gluing data $\alpha^{(\ell)}$ and $\alpha^{(\ell')}$ are linear and the function $\beta$ is quadratic (since the functions $\beta^{(\ell)}$ and $\beta^{(\ell')}$ are linear). This is in contrast to the gluing data $\alpha^{(\ell)}$, $\alpha^{(\ell')}$ and $\beta$ of the initial parameterization which are spline functions of degrees 5, 5 and 6, respectively. One can conclude directly from [7] that the $C^1$-smooth space over this initial configuration is not suitable for isogeometric analysis. This matter of fact will be verified in the following subsection by means of $L^2$ approximation.

*Example 4 – Footprint.* The initial multi-patch domain $\widetilde{\Omega}$, see Fig. 6(a), is inspired by "Yeti's footprint" shown in [21, Fig. 14] and consists of 25 patches $\widetilde{\Omega}^{(\ell)}$, cf. the G+Smo library [15, 23]. Each patch $\widetilde{\Omega}^{(\ell)}$ is parameterized by a geometry mapping $\widetilde{\boldsymbol{F}}^{(\ell)} \in \mathcal{S}_1^{3,2}([0,1]^2) \times$ $\mathcal{S}_1^{3,2}([0,1]^2)$. Fig. 6(b) shows the resulting AS-$G^1$ multi-patch parameterization $\boldsymbol{F}$ by using the space $\mathcal{S}_1^{3,1}([0,1]^2) \times \mathcal{S}_1^{3,1}([0,1]^2)$ for the single parameterizations $\boldsymbol{F}^{(\ell)}$, and Table 1 reports the resulting relative $L^2$ error.

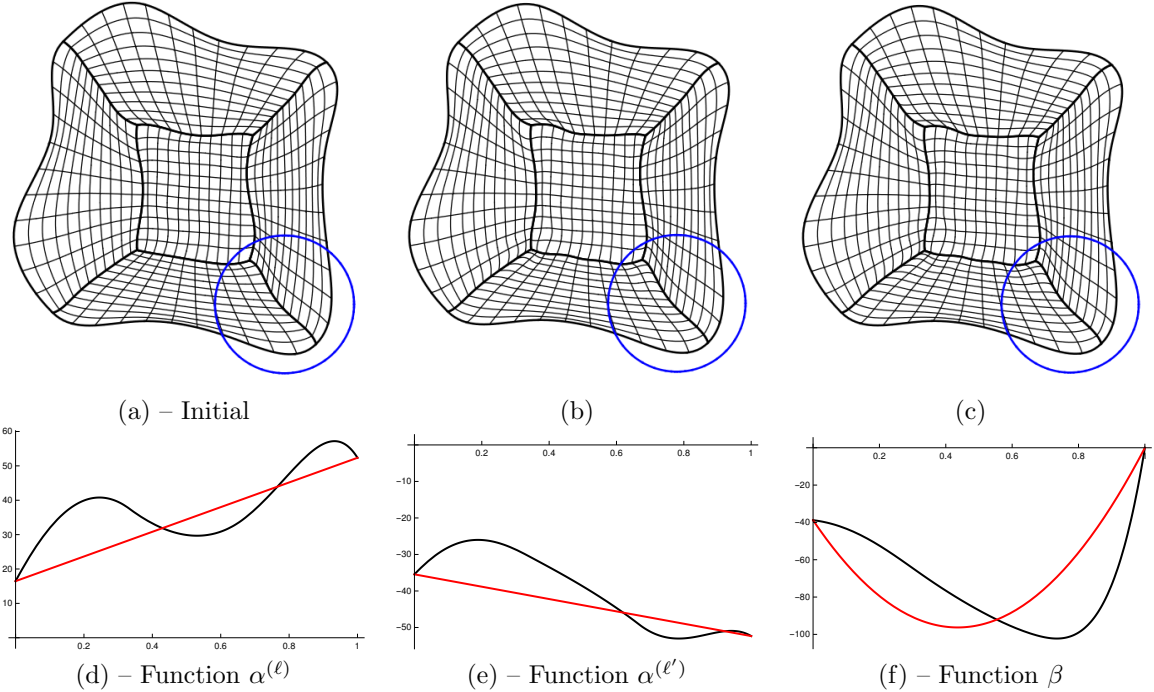

Figure 5: Example 3 – Deformed circle. (a)-(c): Plots of the initial parameterization and of two constructed AS-$G^1$ parameterizations, compare Table 1. (d)-(f): Plots of the gluing data $\alpha^{(\ell)}$, $\alpha^{(\ell')}$ and $\beta$ for the initial parameterization (in black) and for the resulting AS-$G^1$ parameterizations (in red) for the interface in the blue encircled region.

The examples indicate that under mild conditions ($p \geq 3$, $r \leq p-2$ and $k \geq 2$), there always exists a multi-patch parameterization $\boldsymbol{F}$, such that conditions (a)-(c) are satisfied. However, it is unclear whether or not (d) can be satisfied. We conjecture that if $p$ and $k$ are sufficiently large, (a)-(d) always has a solution. In this context, conditions (a), (b) and (d) are natural conditions that are necessary for multi-patch parameterizations that possess optimal approximation properties, cf. [7]. Condition (c) allows us to directly compute the gluing data $\alpha^{(\ell)}$, $\alpha^{(\ell')}$, etc. which linearizes the problem. Consequently this allows us to construct such parameterizations.

### 4.2. $L^2$ approximation

We perform $L^2$ approximation on several multi-patch parameterizations from the previous subsection.

*Refineable $C^1$ isogeometric spaces.* For a given spline space $\mathcal{S}_k^{p,r}([0,1]^2)$, we denote by $h$ the mesh size of the uniform knot mesh in the parameter domain, i.e. $h = \frac{1}{k+1}$. In the following, we write $\mathcal{S}_h^{p,r}([0,1]^2)$ instead of $\mathcal{S}_k^{p,r}([0,1]^2)$. Consider a planar multi-patch domain $\Omega = \cup_{\ell=1}^P \Omega^{(\ell)}$, where each patch $\Omega^{(\ell)}$ is parameterized by a geometry mapping $\boldsymbol{F}^{(\ell)} \in \mathcal{S}_{h_0}^{p,r} \times \mathcal{S}_{h_0}^{p,r}$. For each selection $h = h_0/2^L$, $L \in \mathbb{Z}_0^+$, we have $\boldsymbol{F}^{(\ell)} \in \mathcal{S}_{h}^{p,r} \times \mathcal{S}_{h}^{p,r}$. Here $L$ is the *level of refinement*. We denote by $V_h^1$ the space of $C^1$-smooth isogeometric functions on $\Omega$, which is given by

$$\mathcal{V}_h^1 = \{\phi \in C^1(\Omega) \mid \phi \circ \boldsymbol{F}^{(\ell)} \in \mathcal{S}_h^{p,r}([0,1]^2),\ \ell = 1, \ldots, P\}.$$

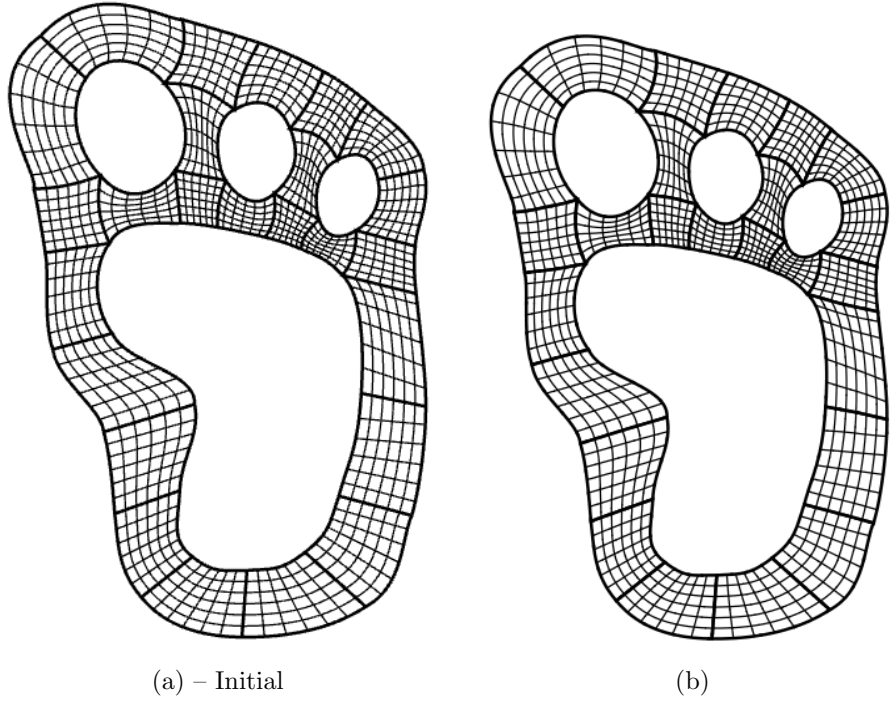

(a) – Initial                                  (b)

Figure 6: Example 4 – Footprint. Compare Table 1.

|     | $p$ | $r$ | $k$ | $\lambda_L$ | $\lambda_U$ | rel $L_2$-err. | $\mathcal{F}_2$ | $\mathcal{F}_L$ | $\mathcal{F}_U$ |
|-----|-----|-----|-----|-------------|-------------|----------------|-----------------|-----------------|-----------------|
| Example 1 – Part of car ||||||||||
| (a) | 2 | $\infty$ | 0 | - | - | - | - | 165.95 | 231.08 |
| (b) | 3 | 1 | 1 | 0 | 0 | $1.03\ 10^{-5}$ | $6.08\ 10^{-4}$ | 165.68 | 256.29 |
| Example 2 – Puzzle piece ||||||||||
| (a) | 3 | 1 | 2 | - | - | - | - | 315.87 | 641.17 |
| (b) | 3 | 1 | 2 | $10^{-3}$ | $5\ 10^{-4}$ | $9.26\ 10^{-6}$ | $5.6\ 10^{-3}$ | 314.54 | 614.2 |
| Example 3 – Deformed circle ||||||||||
| (a) | 3 | 2 | 2 | - | - | - | - | 547.52 | 2422.14 |
| (b) | 3 | 1 | 5 | $10^{-5}$ | $10^{-6}$ | $1.66\ 10^{-5}$ | $2.42\ 10^{-3}$ | 549.17 | 5355.62 |
| (c) | 4 | 1 | 2 | $10^{-5}$ | $10^{-6}$ | $1.75\ 10^{-5}$ | $2.55\ 10^{-3}$ | 549.12 | 5005.05 |
| Example 4 – Footprint ||||||||||
| (a) | 3 | 2 | 1 | - | - | - | - | 58.59 | 115.45 |
| (b) | 3 | 1 | 1 | $10^{-3}$ | $10^{-3}$ | $4.31\ 10^{-6}$ | $1.97\ 10^{-3}$ | 58.19 | 103.27 |

Table 1: Examples 1–4: The selected parameters $p$, $r$ and $k$ for the spline space $\mathcal{S}_k^{p,r}$, the choices of the weights $\lambda_L$ and $\lambda_U$, and the resulting $L^2$ error such as the value of the functional $\mathcal{F}_2$, the value of the parametric length functional $\mathcal{F}_L$ and the value of the uniformity functional $\mathcal{F}_U$ in the objective function (15).

By selecting the sequence $(h_0, h_0/2, h_0/4, \ldots, h_0/2^L, \ldots)$ for the mesh size $h$, we are able to generate a sequence of nested spaces $\mathcal{V}_h^1$, i.e.

$$\mathcal{V}_{h_0}^1 \subset \mathcal{V}_{h_0/2}^1 \subset \mathcal{V}_{h_0/4}^1 \subset \cdots \subset \mathcal{V}_{h_0/2^L}^1 \subset \cdots.$$

A general framework to construct a basis of the $C^1$ isogeometric space $\mathcal{V}_h^1$ was described in [20, Section 2.3]. There it was shown that the construction of a basis is equivalent to finding a null space of a matrix from a homogeneous linear system, which is derived from the required $G^1$ constraints of the graph surface of an isogeometric function along the common interfaces. To get well conditioned bases for the spaces $\mathcal{V}_h^1$ for our examples, we have used the method proposed in [19, Section 6.1] and [18, Section 5.1]. There it was introduced to find bases for $C^2$-smooth isogeometric functions on bilinear multi-patch domains. The method is based on the concept of finding minimal determining sets (cf. [22, Section 5.6]) for the unknown coefficients of the homogeneous linear system and can be used in an analogous way for the space $\mathcal{V}_h^1$.

*Model problem.* Let $N_h = \dim \mathcal{V}_h^1$ and let the set of functions $\{\phi_i\}_{i=0}^{N_h-1}$ be a basis of the space $\mathcal{V}_h^1$. We approximate the function

$$z : \Omega \to \mathbb{R}, \quad z(\boldsymbol{x}) = z(x_1, x_2) = 2\cos(x_1)\sin(x_2), \tag{23}$$

see Fig. 7(top row), by the function

$$u_h(\boldsymbol{x}) = \sum_{i=0}^{N_h-1} c_i \phi_i(\boldsymbol{x}), \quad c_i \in \mathbb{R},$$

by means of the least-squares approach. This means that we minimize the following objective function

$$||u_h - z||_{L^2}^2 = \int_\Omega (u_h(\boldsymbol{x}) - z(\boldsymbol{x}))^2 \mathrm{d}\boldsymbol{x} \to \min_{c_i},$$

This minimization problem is equivalent to solving the system of linear equations

$$M\boldsymbol{c} = \boldsymbol{z}, \quad \boldsymbol{c} = (c_i)_{i=0}^{N_h-1},$$

where the entries of the matrix $M = (m_{i,j})_{i,j=0}^{N_h-1}$ and of the vector $\boldsymbol{z} = (z_i)_{i=0}^{N_h-1}$ are given by

$$m_{i,j} = \int_\Omega \phi_i(\boldsymbol{x})\phi_j(\boldsymbol{x}) \mathrm{d}\boldsymbol{x}, \text{ and } z_i = \int_\Omega z(\boldsymbol{x})\phi_i(\boldsymbol{x}) \mathrm{d}\boldsymbol{x},$$

respectively. The isogeometric approach allows to rewrite the entries $m_{i,j}$ and $z_i$ as

$$m_{i,j} = \sum_{\ell=1}^{P} \int_{[0,1]^2} (\phi_i \circ \boldsymbol{F}^{(\ell)}(u,v))(\phi_j \circ \boldsymbol{F}^{(\ell)}(u,v))|\det(J^{(\ell)}(u,v))|\mathrm{d}u\mathrm{d}v,$$

and

$$z_i = \sum_{\ell=1}^{P} \int_{[0,1]^2} z(\boldsymbol{F}^{(\ell)}(u,v))(\phi_i \circ \boldsymbol{F}^{(\ell)}(u,v))|\det(J^{(\ell)}(u,v))|\mathrm{d}u\mathrm{d}v,$$

respectively, where $J^{(\ell)}$, $\ell \in \{1,\ldots,P\}$, is the Jacobian of $\boldsymbol{F}^{(\ell)}$.

*Numerical examples.* We perform $L^2$ projection on the isogeometric space defined from the three AS-$G^1$ multi-patch parameterizations, shown in Figs. 3(b)-5(b), and on their associated initial parameterizations (see Figs. 3(a)-5(a)). To approximate the function (23), see Fig. 7 (top row), we have used nested spline spaces $\mathcal{V}_h^1$ of bidegree $(3,3)$ for the different levels $L$ of refinement. The resulting relative $L^2$ (i.e. $H^0$), $H^1$ and $H^2$ errors are reported in Fig. 7 (second to fourth row), and the resulting absolute errors are plotted in Figs. 8-10.

In case of the AS-$G^1$ multi-patch parameterizations, the estimated convergence rates, see Fig. 7 (second to fourth row), are of order $\mathcal{O}(h^{4-i})$ in the $H^i$-norms, $i = 0, 1, 2$. These optimal convergence rates and the resulting absolute errors, see Figs. 8-10, verify the optimal approximation properties of the $C^1$ isogeometric spaces over AS-$G^1$ multi-patch parameterizations.

In case of the initial (non-AS-$G^1$) multi-patch parameterization, the numerical results indicate for all three examples no convergence at all in $L^\infty$ along most interfaces and consequently the minimal convergence order $\mathcal{O}(h^{1/2})$ in $L^2$ towards the exact solution. This can be observed by the estimated convergence rates, compare Fig. 7 (second to fourth row), and by the resulting absolute errors, see Fig. 8-10, which do not decrease along some patch interfaces. The $H^1$ and $H^2$ errors grow with mesh refinement, which is to be expected since we are approximating by $L^2$-projection in all cases. The poor convergence behaviour for non-AS-$G^1$ multi-patch parameterizations derives from the fact that the $C^1$ isogeometric spaces are in general restricted to low degree ($\ll p$) splines or polynomials along the patch interfaces, cf. [7]. In the examples presented here, only linear functions remain along the interfaces.

## 5. Conclusion

We have presented a method to construct a planar multi-patch parameterization, which belongs to the class of so-called analysis suitable $G^1$ (in short, AS-$G^1$) multi-patch parameterizations. This class of parameterizations contains exactly those multi-patch parameterizations, which allow to define $C^1$ isogeometric spaces with optimal approximation properties (cf. [7]). AS-$G^1$ multi-patch parameterizations have to satisfy specific constraints, called AS-$G^1$ constraints, which are highly non-linear conditions for unknown geometries. We linearize these constraints and introduce additional ones to construct an AS-$G^1$ multi-patch parameterization, which is as close as possible to an initial planar multi-patch parameterization, and which additionally interpolates the boundary, the vertices and the first derivatives at the vertices of this initial parameterization.

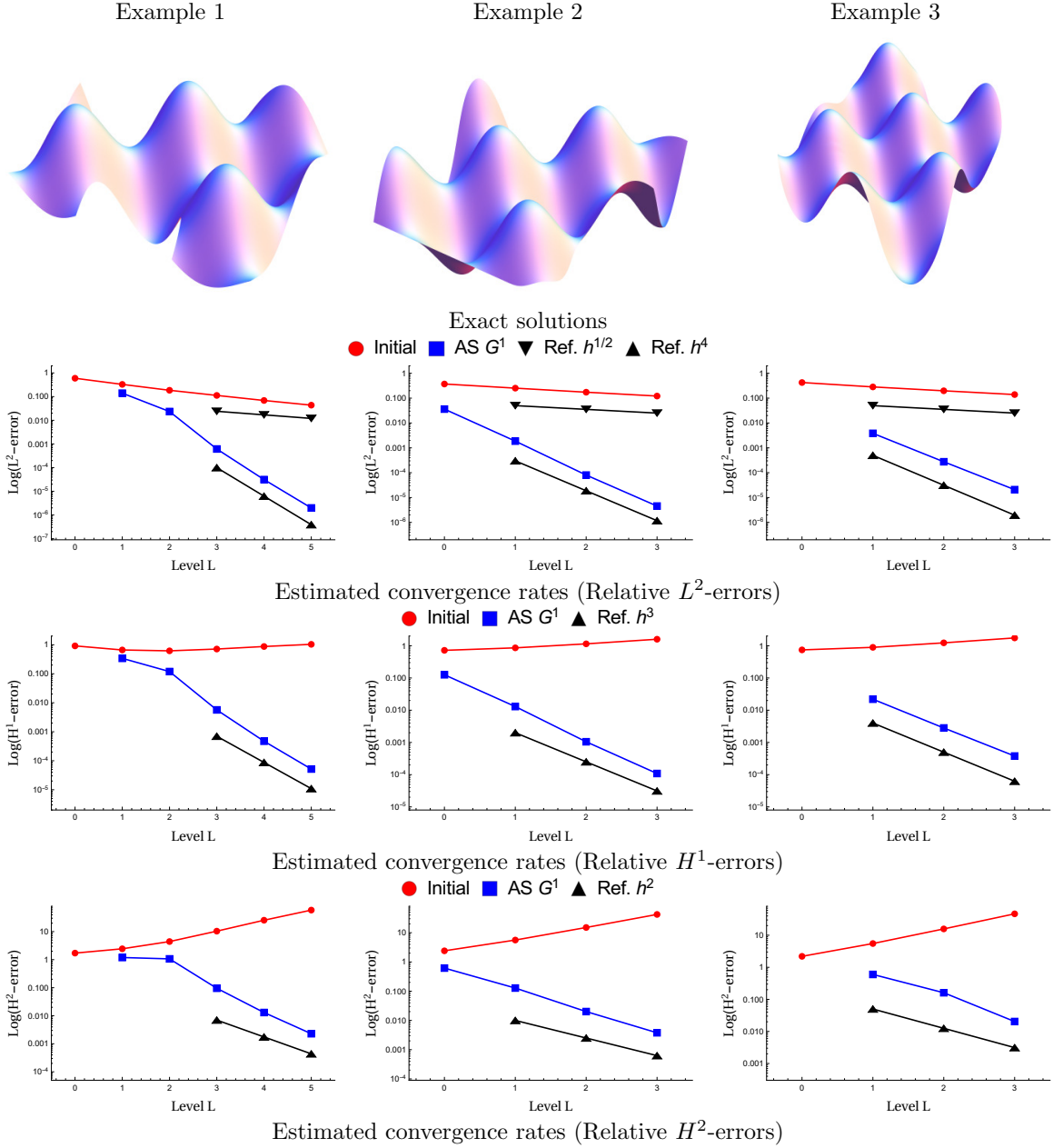

Figure 7: Example 1–3: Performing $L^2$ approximation on the initial and AS-$G^1$ multi-patch parameterizations. Compare Figs. 8–10.

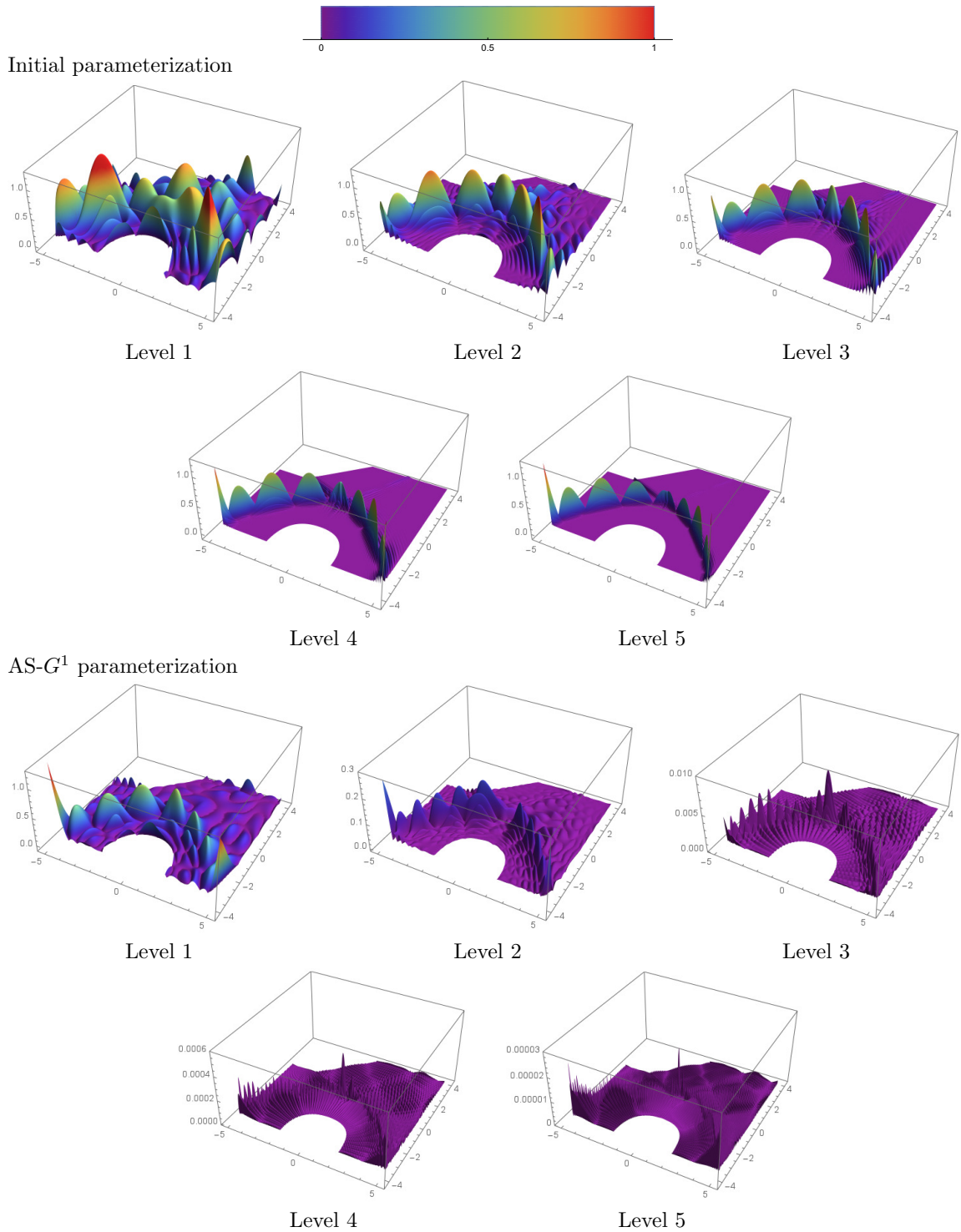

Figure 8: $L^2$ approximation – Part of car (Example 1): Absolute errors with respect to the exact solution (23) (see Fig. 7 (top row)) for different levels $L$ by performing $L^2$ approximation on the initial multi-patch parameterization (see Fig. 3(a)) and on the constructed AS-$G^1$ multi-patch parameterization (see Fig. 3(b)).

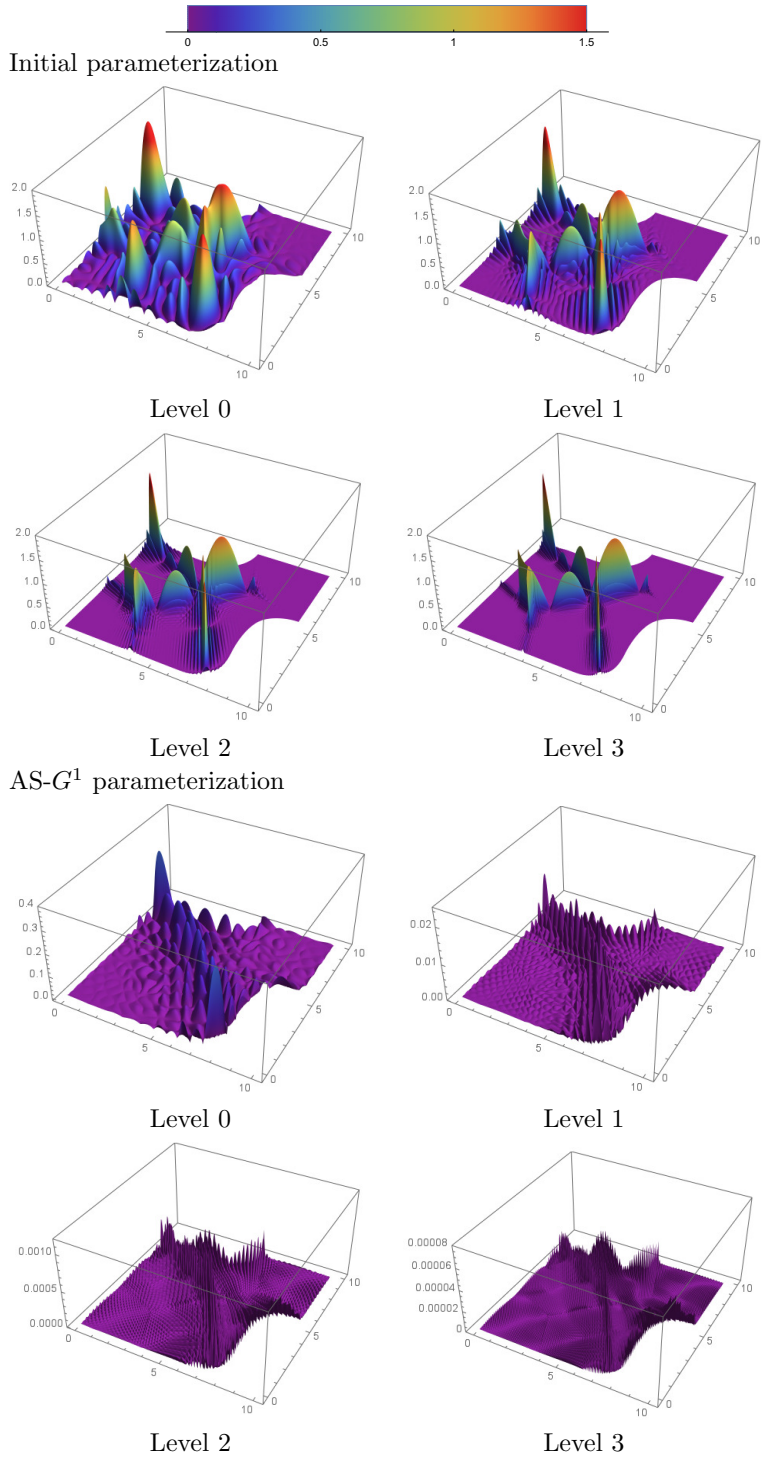

Figure 9: $L^2$ approximation – Puzzle piece (Example 2): Absolute errors with respect to the exact solution (23) (see Fig. 7 (top row)) for different levels by performing $L^2$ approximation on the initial multi-patch parameterization (see Fig. 4(a)) and on the constructed AS-$G^1$ multi-patch parameterization (see Fig. 4(b)).

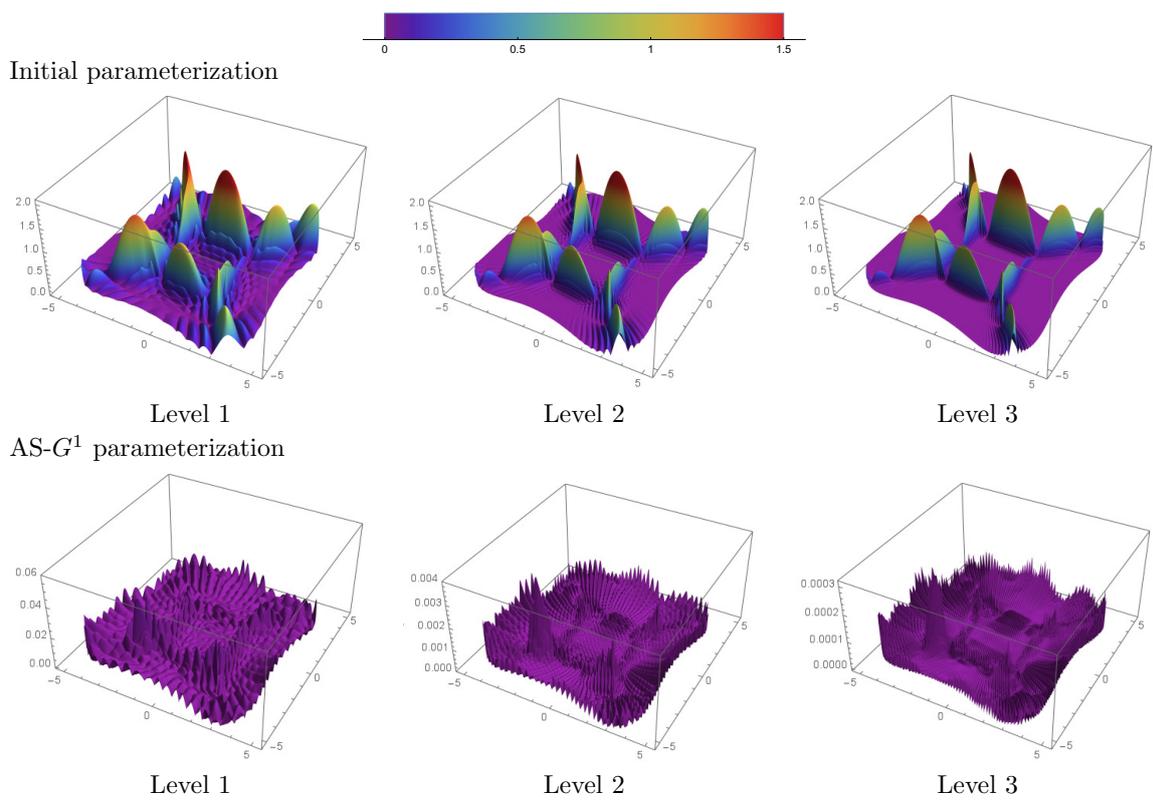

Figure 10: $L^2$ approximation – Deformed circle (Example 3): Absolute errors with respect to the exact solution (23) (see Fig. 7 (top row)) for different levels by performing $L^2$ approximation on the initial multi-patch parameterization (see Fig. 5(a)) and on the constructed AS-$G^1$ multi-patch parameterization (see Fig. 5(b)).

Our proposed method is simple and requires only to solve a system of linear equations. We have presented several examples where we constructed AS-$G^1$ planar multi-patch parameterizations. On the one hand we could demonstrate the potential of our algorithm to construct such reparameterizations reproducing any planar multi-patch geometry exactly. On the other hand we could show that the AS-$G^1$ constraints themselves and also the linearization of these constraints are not a real restriction to generate a multi-patch parameterization of good quality from a given multi-patch domain, which proves the flexibility of AS-$G^1$ geometries.

Since the AS-$G^1$ constraints are highly sensitive with respect to small perturbations, all computations have been performed symbolically. Our current and further research will focus on the reduction of this sensitivity by generating a well conditioned basis, containing elements with small local supports, for the solution space of the AS-$G^1$ constraints. This is a challenging problem and an essential ingredient for the construction of a well-conditioned and locally supported basis for the $C^1$ isogeometric space over AS-$G^1$ multi-patch parameterizations.

Our work has two important implications:

- We give numerical evidence that the AS-$G^1$ definition leaves enough flexibility to allow the actual construction of AS-$G^1$ parameterizations of planar domains, and

- we further confirm the claims of [7], that is, optimal convergence is precluded in general on arbitrary parameterizations (no convergence in $L^\infty$ is evidenced in some cases) but it is guaranteed on AS-$G^1$ parameterizations.

Further topics which are worth to study are the use of such parameterizations in more complex isogeometric simulations (e.g., for shells) and the extension of our algorithm to multi-patch surfaces or volumetric multi-patch domains.

**Acknowledgment.** The first two authors (M. Kapl and G. Sangalli) were partially supported by the European Research Council through the FP7 ERC Consolidator Grant n.616563 HIGEOM, and by the Italian MIUR through the PRIN Metodologie innovative nella modellistica differenziale numerica. This support is gratefully acknowledged.